\documentclass[12pt]{elsart}
\usepackage{amssymb}
\newtheorem{theo}{Theorem}
\newcommand{\dt}{{\partial\over \partial t}}
\newcommand{\dx}{{\partial\over \partial x}}
\newcommand{\du}{{\partial\over \partial u}}
\begin{document}

\begin{frontmatter}
\title{Euler-Bernoulli beams from a symmetry standpoint-characterization
   of equivalent equations}

\author{C\'elestin Wafo Soh\thanksref{w}}

\address{Mathematics Department,
 College of Science, Engineering, and Technology  Jackson State
 University,
 JSU Box 17610, 1400 J R Lynch St., \\Jackson, MS 39217, USA}
\ead{wafosoh@yahoo.com, celestin.wafo@jsums.edu}
 \thanks[w]{To the loving  memory of my brother L\'eopold Fotso Simo.}

\begin{abstract}
We completely solve the equivalence problem for Euler-Bernoulli
equation using Lie symmetry analysis. We show that the quotient of
the symmetry Lie algebra of the Bernoulli equation by the
infinite-dimensional Lie algebra spanned by solution symmetries is a
representation of one of the following Lie algebras: $2A_1$,
$A_1\oplus A_2$, $3A_1$, or $A_{3,3}\oplus A_1$. Each quotient
symmetry Lie algebra determines an equivalence class of
Euler-Bernoulli equations. Save for the generic case corresponding
to arbitrary lineal mass density and flexural rigidity, we
characterize the elements of each class by giving a determined set
of differential equations satisfied by physical parameters (lineal
mass density  and flexural rigidity). For each class,  we provide a
simple representative and we explicitly construct transformations
that maps a class member to its representative. The maximally
symmetric class described by the four-dimensional quotient symmetry
Lie algebra $A_{3,3}\oplus A_1$ corresponds to Euler-Bernoulli
equations homeomorphic to the uniform one (constant lineal mass
density  and flexural rigidity) . We rigorously derive some
non-trivial and non-uniform Euler-Bernoulli equations reducible to
the uniform unit beam. Our models extend and emphasize the symmetry
flavor of Gottlieb's iso-spectral beams (Proceedings of the Royal
Society London A 413 (1987) 235-250)
\end{abstract}

\end{frontmatter}
\section{Introduction}
Da Vinci and Galileo foresaw the need for a theory of vibrating thin
beams. However they suggested theories that were either incomplete
or erroneous. Da Vinci's theory was more descriptive and based on
detailed sketches rather than physical laws and equations : he
lacked tools such as Hooke's law, Newton's laws, and calculus which
postdate him. In Galileo's approach, the nemesis was an incorrect
calculation of the load carrying capacity of transversely loaded
beams. We owe the first consistent thin beams theory to the
Bernoullis. Jacob Bernoulli developed an elasticity theory in which
the curvature of an elastic beam is proportional to its bending
moment. Relying on his uncle elasticity theory, Daniel Bernoulli
derived a partial differential equation governing the motion of a
thin vibrating beam. Leonard Euler extended and applied the
Benoullis theory to loaded beams.

In Euler-Bernoulli beam theory, the transversal motion of an
unloaded thin elastic beam is governed by the partial differential
equation
\begin{equation}
{\partial^2\over \partial \,x^2}\left ( f(x)\, {\partial^2
\,u\over
\partial \,x^2} \right ) +m(x)\,{\partial^2 u\over \partial \,t^2} = 0,
\quad t>0,\quad 0<x<L, \label{i1}
\end{equation}
where $f(x)>0$ is the flexural rigidity, $m(x)>0$ is the lineal mass
density, and $u(t,x)$ is the transversal displacement at time $t$
and position $x$ from one end of the beam taken as origin. Equation
(\ref{i1}) must be solved subject to initial and boundary conditions
such as clamped ends, hinged ends, and  free ends boundary
conditions.

In this paper, our focus is on the equivalence problem for Eq.
(\ref{i1}): we seek necessary and sufficient conditions under which
two equations of the form (\ref{i1}) can be mapped to each other
using an invertible change of the dependent and independent
variables. A particular case of this equivalence problem was tackled
by Gottlieb \cite{gott} who was interested in equations of the form
(\ref{i1}) that are equivalent to the uniform ( constant $f$ and
$m$) beam equation. In the same vain, Bluman and Kumei
\cite{blku}(Chapter 6, Section 6.5) studied the problem of reducing
a linear partial differential equation to a constant coefficient
one.

The layout of this paper is the following. There are four sections
including this introduction. Section 2 deals with the complete Lie
symmetry classification of Eq. (\ref{i1}). Section 3 is dedicated to
the construction of equivalence transformations. We recapitulate our
findings in Section 4.

\section{Symmetry analysis
of Euler-Bernoulli equation}

Our goal in this section is to study the symmetry breaking of Eq.
(\ref{i1}). We assume that the reader is familiar with the rudiments
of Lie's symmetry theory \cite{blku,ovsi,olve}.

A vector field
\begin{equation}
X =\tau (t,x,u)\,\dt +\xi (t,x,u)\,\dx +\eta(t,x,u)\,\du
\label{lie1}
\end{equation}
is a Lie symmetry of Eq. (\ref{i1}) if
\begin{equation}
X^{[4]} \left. \left (  {\partial^2\over \partial \,x^2}\left (
f(x)\, {\partial^2 \,u\over
\partial \,x^2} \right ) +m(x)\,{\partial^2 u\over \partial \,t^2}  \right )\right |_{Eq. (\ref{i1})} =0, \label{lie2}
\end{equation}
where $X^{[4]}$ is the fourth prolongation of $X$ which is
calculated using the formulas
\begin{eqnarray}
X^{[k]} &=& X+\sum_{1\le |J|\le k} \eta_{J}\,{\partial \over
\partial
u_J},\quad  J=(j_1,j_2),\quad |J|=j_1+j_2, \label{lie3'}\\
\quad u_J
&= & \partial^{|J|}u/\partial t^{j_1}\partial x^{j_2}, \label{lie3} \\
\eta_{J}  & =& D_{J} (\eta -\tau u_t-\xi u_x)+\tau\,
u_{J,t}+\xi\,u_{J,x},\quad u_{J,r}=\partial u_J/\partial r,
\label{lie4}\\
 D_J & = &  D_{j_1}\,D_{j_2},\quad D_{j_1} =(D_t)^{j_1},\quad
 D_{j_2}=(D_x)^j_2, \label{lie5} \\
 D_{t} &=& \dt +u_t\du + u_{tt}{\partial \over \partial u_t} + u_{tx}{\partial \over \partial u_x} +
 \cdots \label{lie6}\\
 D_{x} &=& \dx +u_x\du + u_{tx}{\partial \over \partial u_t} + u_{xx}{\partial \over \partial u_x} +
 \cdots \label{lie7}
 \end{eqnarray}

Since the symmetry coefficients $\tau$, $\xi$ and $\tau$ do not
depend explicitly on the derivatives of $u$, the left-hand side of
Eq. (\ref{lie1}) is a polynomial in the derivatives of $u$. Thus we
may set its coefficients to zero to obtain an over-determined system
of linear partial differential equations. In order to avoid the
appearance of the integral $\int (m/f)^{1/4} dx$ in our
calculations, we  express the lineal mass density  as follows.
\begin{equation}
m(x)= (g'(x))^4\,f(x). \label{d0}
\end{equation}

After some calculations, the determining equations for the
symmetries  simplify  to

\begin{eqnarray}
& & \tau   =  4\,c_1\,t+c_2 \label{d1}\\
& & \xi   =  {2\,c_1\,g\over g'}+{2\,c_3\over g'}, \label{d2}\\
& & \eta   = -\left ( {c_1\,f'\,g\over f\,g'}+{c_3\,f'\over
f\,g'}+{3\,c_1\,g\,g''\over g'^2}+{3\,c_3\,g''\over g'^2}+c_4 \right
)u + a(t,x)
\label{d3} \\
 & & ( f a_{xx})_{xx} + g'^4 f a_{tt} =0 \label{d3'}\\
 & & c_1 H_{11}+c_3H_{12} =0 ,\label{d4}\\
 & & c_1 H_{21}+c_3H_{22} =0, \label{d5}\\
 & &  (c_1 H_{11}+c_3H_{12})_{xx} =0, \label{d6}
 \end{eqnarray}
where $c_1$ to $c_4$ are integration constants, and the differential
functions $H_{ij}$ are relegated to Appendix \ref{ap1}  due to their
size. It can be readily seen that Eq. (\ref{d6}) is a mere
differential consequence of Eq. (\ref{d4}). Thus the equations we
have to solve are Eqs. (\ref{d4}) and (\ref{d5}).

It is well-known (see for example \cite{blku}) that the symmetry Lie
algebra of a scalar linear partial differential equation is of the
form $L^r\oplus L^{\infty}$, where $L^r$ is a finite-dimensional Lie
algebra and $L^{\infty}$ is an infinite-dimensional ideal of the
symmetry Lie algebra spanned by the so-called solution symmetries.
In our case
$$L^{\infty}\;=\;<a(t,x)\partial_u >,$$ where $a(t,x)$ solves Eq.
(\ref{d3'}) i.e.  Euler-Bernoulli equation. Thus we are left with
characterizing $L^r$. The determining equations contain four
constants {\it viz.} $c_1$ to $c_4$ with two constants that are
unconstrained. Therefore the dimension of $L^r$, $r$, lies between
two and four. Below we elucidate all the possibilities.

\subsection*{Case I: $c_1$ and  $c_3$ are arbitrary constants}
In this case, Eqs. (\ref{d4})-(\ref{d5}) split into the following
system.

\begin{eqnarray}
H_{11} &= & 0, \label{d7}\\
H_{12} &= & 0, \label{d8}\\
H_{21} &= & 0, \label{d9}\\
H_{22} &= & 0. \label{d10}
\end{eqnarray}
The system (\ref{d7})-(\ref{d10}) is an over-determined system of
nonlinear ordinary differential equations for $f$ and $g$. Thus,
{\it a priori}, we are not guaranteed of a solution. Replace Eq.
(\ref{d7}) by the combination Eq.(\ref{d7})-$g
\times$Eq.(\ref{d8}). Solve the resulting equation for $g^{(3)}$
to obtain
\begin{equation}
g^{(3)}= {3\over 10}\, {g'\,f'^2\over f^2}-{2\over
5}\,{g'\,f''\over f}+{3\over 2}\,{g''^2\over g'}\;\cdot
\label{d11}
\end{equation}
Use Eq. (\ref{d11}) to eliminate the derivative $g^{(3)}$ and
$g^{(4)}$ from Eq. (\ref{d8}). It results the equation
\begin{equation}
f^{(4)}={f'\,f^{(3)}\over f}+{11\over 10}\,{f''^2\over f}-{12\over
5}\,{f'^2\,f''\over f^2}+{9\over 10}\,{f'^4\over f^3}\,\cdot
\label{d12}
\end{equation}
Employing Eqs (\ref{d11})-(\ref{d12}) to get rid of  the derivatives
$g^{(3)}$ to $g^{(6)}$, $f^{(4)}$, and $f^{(5)}$ from Eq.
(\ref{d9})-(\ref{d10}), we discover that Eqs.(\ref{d9})-(\ref{d10})
are identically satisfied. To sum up, we have established that the
over-determined system of Eqs.(\ref{d7})-(\ref{d10}) is equivalent
to the determined system formed by Eqs. (\ref{d11})-(\ref{d12}).
Thus provided Eqs. (\ref{d11})-(\ref{d12}) are satisfied, the finite
part of the symmetry Lie algebra, $L^4$, is spanned by the operators
\begin{eqnarray}
 X_1 & =& \dt,\quad X_2 =u\;\du,\label{d13}\\
 X_3 & = & 4t\,\dt+ {2\,g\over g'}\,\dx-\left( {f'\,g\over fg'}+{3\,g\,g''\over g'^2} \right )u\,\du
 \label{d14}\\
 X_4 & = & {2\over g'}\,\dx -\left ( {f'\over f\,g'}+{3\,g''\over g'^2}  \right )u\,\du
 \label{d15}
\end{eqnarray}
Simple computations show that the nonzero commutators of the
symmetry generators are
$$ [X_1,X_3]=4X_1,\quad [X_3,X_4]=-2X_4. $$
By making the change of basis
$$ e_1=X_1,\quad e_2 ={1\over \sqrt{2}}\;X_4,\quad e_3={1\over 4}\;X_3,
\quad e_4=X_2,$$ it can be seen that $L^4$ is equivalent to
$A_{3,3}\oplus A_1$ in Patera and Winternitz \cite{pawi}
classification scheme.
\subsection*{Case II: $ c_1$ is arbitrary and $c_3=0$}
 The determining equations (\ref{d4})-(\ref{d5}) become
 \begin{eqnarray}
 H_{11} &= & 0, \label{d16}\\
 H_{21} & = & 0. \label{d17}
 \end{eqnarray}
 We aim at rewriting the constraints (\ref{d16})-(\ref{d17}) in
 terms of lowest possible derivatives of $f$ and $g$. In order to
 achieve this goal, we first solve Eq. (\ref{d16}) for $g^{(4)}$ to
 obtain
 \begin{eqnarray}
 g^{(4)} &=& 6\,{g''\,g^{(3)}\over g'}-2\,{g'\,g^{(3)}\over g}-{2\over
 5}\,{g'\,f^{(3)}\over f} +6\,{g''^3\over g'^2} +3\,{g''^2\over g}
 +{4\over 5}\,{g''\,f''\over f}\nonumber \\
 & & -{3\over 5}\,{f'^2\,g''\over g}+{g'\,f'\,f''\over f}-{4\over  5}\,{g'^2\,f''\over g\,f}
 -{3\over 5}\,{g'\,f'^3\over f^3}+{3\over 5}\,{g'^2\,f'^2\over f^2\,g}\; \cdot \label{d18}
 \end{eqnarray}
Employ  Eq. (\ref{d18}) to get rid of the derivatives $g^{(4)}$ to
$g^{(6)}$ from Eq. (\ref{d17}) and solve the resulting equation for
$f^{(5)}$ to obtain
\begin{eqnarray}
f^{(5)}&=& \frac{-18\,{f'}^5}{5\,{f}^4} +
  \frac{18\,{f'}^4\,g'}{5\,{f}^3\,g} +
  \frac{54\,{f'}^3\,f''}{5\,{f}^3} -
  \frac{48\,{f'}^2\,g'\,f''}{5\,{f}^2\,g} -
  \frac{7\,f'\,{f''}^2}{{f}^2} +
  \frac{22\,g'\,{f''}^2}{5\,f\,g} \nonumber \\
  &&-
  \frac{18\,{f'}^4\,g''}{5\,{f}^3\,g'}
   + \frac{48\,{f'}^2\,f''\,g''}{5\,{f}^2\,g'} -
  \frac{22\,{f''}^2\,g''}{5\,f\,g'} -
  \frac{22\,{f'}^2\,f^{(3)}}{5\,{f}^2} +
  \frac{4\,f'\,g'\,f^{(3)}}{f\,g} \nonumber \\
  & & +
  \frac{16\,f''\,f^{(3)}}{5\,f}
    -
  \frac{4\,f'\,g''\,f^{(3)}}{f\,g'} +
  \frac{2\,f'\,f^{(4)}}{f} - \frac{4\,g'\,f^{(4)}}{g} +
  \frac{4\,g''\,f^{(4)}}{g'} \; \cdot
 \label{d19}
\end{eqnarray}
  Thus the constraints Eqs. (\ref{d16})-(\ref{d17}) are equivalent to
 the simplified constraints (\ref{d18})-(\ref{d19}). Provided $f$
 and $g$ fulfilled  Eqs.(\ref{d18})-(\ref{d19}), the finite part of the
 symmetry Lie algebra, $L^{3,1}$, is spanned by $X_1$, $X_2$, and
 $X_3$. This Lie algebra corresponds to $A_1\oplus A_2$ in Patera
 and Winternitz \cite{pawi} classification of lower-dimensional Lie
 algebras.
\subsection*{Case III: $c_1=0$ and $c_3$ is arbitrary}
Here, Eqs. (\ref{d16})-(\ref{d17}) are equivalent to the system of
equations
\begin{eqnarray}
H_{21} & = & 0, \label{d20}\\
H_{22} & = & 0. \label{d21}
\end{eqnarray}
By following the same {\it modus operandi} as in Case II, we arrive
at the following simplified constraints on $f$ and $g$.

\begin{eqnarray}
g^{(4)} &= & \frac{-3\,{f'}^3\,g'}{5\,{f}^3} +
\frac{f'\,g'\,f''}{{f}^2} - \frac{3\,{f'}^2\,g''}{5\,{f}^2} +
\frac{4\,f''\,g''}{5\,f} \nonumber  \\
& & -
  \frac{6\,{g''}^3}{{g'}^2} - \frac{2\,g'\,f^{(3)}}{5\,f}
+ \frac{6\,g''\,g^{(3)}}{g'}\, ,\label{d22}\\
f^{(5)} &=&\frac{-18\,{f' }^5}{5\,{f  }^4} +
\frac{54\,{f'}^3\,f''}{5\,{f}^3} - \frac{7\,f'\,{f''}^2}{{f}^2} -
\frac{18\,{f'}^4\,g''}{5\,{f}^3\,g'} \nonumber \\
& &+
  \frac{48\,{f'}^2\,f''\,g''}{5\,{f}^2\,g'}
  - \frac{22\,{f''}^2\,g''}{5\,f\,g'}
  - \frac{22\,{f'}^2\,f^{(3)}}{5\,{f}^2} \nonumber \\
  & & +
  \frac{16\,f''\,f^{(3)}}{5\,f}
  - \frac{4\,f'\,g''\,f^{(3)}}{f\,g'}
  + \frac{2\,f'\,f^{(4)}}{f}
  + \frac{4\,g''\,f^{(4)}}{g'}\, \cdot \label{d23}
\end{eqnarray}
The finite portion of the symmetry Lie algebra, $L^{3,2}$, is
spanned by $X_1$,$X_2$ and $X_4$. The Lie algebra $L^{3,2}$ is
nothing but the three-dimensional Abelian Lie algebra denoted by
$3A_1$ in Patera and Winternitz \cite{pawi} classification of
lower-dimensional Lie algebras.
\subsection*{Case IV: $c_1=0=c_3$ }
This is the generic case: there are no constraints on $f$ and $g$.
The finite part of the Lie algebra, $L^2$, is generated by $X_1$
and $X_2$. It is the two-dimensional Abelian Lie algebra $2A_1$.

To sum up, we have established the following result.
\begin{theo}
Denote the symmetry Lie algebra of Euler-Bernoulli equation by
$S$, and $L^{\infty}$ the infinite-dimensional Lie algebra
generated by the solution symmetries. Then,
$S/L^{\infty}=A_{3,3}\oplus A_1,\,A_1\oplus A_2,\; 3A_1,\;2A_1,$
depending on wether $f$ and $g$ respectively satisfy
Eqs.(\ref{d11})-(\ref{d12}), Eqs. (\ref{d18})-(\ref{d19}), Eqs.
(\ref{d22})-(\ref{d23}), or are arbitrary.
\end{theo}
\section{Equivalence classes and mapping to canonical elements}
From a symmetry standpoint there are essentially four classes of
Euler-Bernoulli equations. Equations of the same class share the
same symmetry structure, and are homeomorphic. Our aim in this
section is to select a simple representative for the non-generic
classes (i.e. all the cases save the case where $f$ and $g$ are
arbitrary), and constructively show how class members are mapped
to their representative.

In the sequel, we shall use capitalized variables to describe
representative of each class. The rational behind this choice shall
be apparent as this section unfold. We denote by $[\mathcal{L}]$ the
set of all Euler-Bernoulli equations having $\mathcal{L}$ as the
finite part of their symmetry Lie algebra. Notations not introduce
here are those of the previous sections.

\subsection{The class $[A_{3,3}\oplus A_1]$}
\subsubsection{Construction of similarity transformations}
It can be readily verified that $f(x)=1$ and $g(x)=x$ satisfy Eqs.
(\ref{d11})-(\ref{d12}). Thus, we select as representative element
of this class the equation
\begin{equation}
U_{XXXX}+U_{TT}=0. \label{e1}
\end{equation}
The finite part of the symmetry Lie algebra of Eq. (\ref{e1}) is
generated by the vectors
\begin{equation}
Y_1 ={\partial \over \partial T},\quad Y_2 =U\,{\partial\over
\partial U},\quad Y_{3,1} =4T\,{\partial\over \partial T}+2X\,{\partial
\over \partial X},\quad Y_4=2\,{\partial \over \partial X}\,\cdot
\label{e1'}
\end{equation}
The invertible transformation
\begin{equation}
T=T(t,x,u),\quad X=X(t,x,u),\quad U=U(t,x,u) \label{e2}
\end{equation}
maps an element of $[A_{3,3}\oplus A_1]$ to Eq. (\ref{e1}) if and
only if the same transformation maps $<X_1,X_2,X_3,X_4>$ to
$<Y_1,Y_2,Y_{3,1},Y_4>$. We look for a transformation (\ref{e2})
such that
\begin{equation}
 X_1\mapsto  Y_1,\quad X_2\mapsto Y_2,\quad
Y_{3,1}\mapsto X_3+\mu_1X_1+\mu_2X_4,\quad Y_4 \mapsto
X_4.\label{e3}
\end{equation}
Recall that the transformation defined by Eq. (\ref{e2}) maps a
vector field $\Gamma$ depending on $t,\,x$ and $u$ to the vector
field $\Gamma(T)\,\partial_T
+\Gamma(X)\,\partial_X+\Gamma(U)\,\partial_U$. Thus, in order to
realize Eqs.(\ref{e3}a)-(\ref{e3}b), and Eq. (\ref{e3}c), we have to
impose
\begin{eqnarray}
X_1(T)& =& 1,\quad X_1(X)=0,\quad X_1(U)=0, \label{e4}\\
X_2(T) &=& 0,\quad X_2(X)=0,\quad X_2(U)=U \label{e5}\\
X_4(T) &=&0,\quad X_4(X)=2,\quad X_4(U)=0\label{e6}
\end{eqnarray}
Solving  Eqs. (\ref{e4})-(\ref{e6}), we obtain
\begin{equation}
T= t+k_1,\quad X=g+k_2,\quad U=k_3\,u \,\sqrt{f\,g'^3}, \label{e7}
\end{equation}
where $k_1$,$k_2$, and $k_3\ne 0$ are constants. It can be readily
verify that Eq. (\ref{e3}b) is satisfied for $\mu_1=k_1$ and
$\mu_2=k_2$. We summarize the findings of this subsection as
follows.
\begin{theo} Equations of $[A_{3,3}\oplus A_1]$ are homeomorphic
to Eq. (\ref{e1}). A transformation that maps an arbitrary element
of $[A_{3,3}\oplus A_1]$ to Eq. (\ref{e1}) is given by Eq.
(\ref{e7}).
\end{theo}
\subsubsection{Examples of non-uniform beams homeomorphic to uniform
beams: generalization and symmetry justification of Gottlieb's
iso-spectral models }

Here we look for closed-form solutions of the uncoupled system of
Eqs. (\ref{d11})-({\ref{d12}).

We may rewrite Eq. (\ref{d11}) as
\begin{equation}
\{g,x\} ={3\over 10}\,{f'^2\over f^2}-{2\over 5}\;{f''\over f},
\label{s1}
\end{equation}
where $ \{y,x\}=y''/y'-(3/2)\,(y''/y')^2$ is the so-called
Schwartzian `derivative' (it is not really a derivative  but a
differential invariant!) of $y$ with respect to $x$. From the
well-known result on inversion of the Schwartzian derivative
\cite{inc}, we infer that
\begin{equation}
g ={y_2(x)\over y_1(x)}\, , \label{s2}
\end{equation}
where $y_1$ and $y_2$ are two linearly independent solutions of the
second-order linear ordinary differential equation
\begin{equation}
y'' +{1\over 20}\;  \left (3\,{f'^2\over f^2}-4\;{f''\over f}\right
)\; y =0. \label{s3}
\end{equation}
We are now left with solving Eq. (\ref{d12}). Its symmetry Lie
algebra is spanned by the operators
\begin{equation}
\Gamma_1 =\dx,\quad \Gamma_2=x\,\dx,\quad \Gamma_2 =f\,{\partial
\over \partial f}\;\cdot \label{s3}
\end{equation}
Since the Lie algebra $<\Gamma_1,\Gamma_2,\Gamma_3>$ is solvable
(its second derived algebra is trivial), we may use successive
reduction to depress the order of Eq. (\ref{d12}) by three to obtain
an Abelian equation  of the second kind which we could not solve in
closed-form. For details about the successive reduction of Eq.
(\ref{d12}), we refer the reader to Appendix \ref{ap2}. Due to the
lack of closed-formed solution of the reduced equation, we look for
an invariant solution of Eq. (\ref{d12}). The most general invariant
solution under a linear combination of $\Gamma_1$ to $\Gamma_3$ are
\begin{equation}
f =K(Ax+B)^m, \quad f=C\e^{Dx} \label{s4}
\end{equation}
where $K$,  $A$ to $D$, and $m$ are constants.  Note that the second
ansatz Eq. (\ref{s4}b) satisfies Eq. (\ref{d12}) if and only if
$D=0$, and this possibility is already included in the first ansatz.
The substitution of Eq. (\ref{s4}a) into Eq. (\ref{d12}) yields the
constraint
\begin{equation}
m(4\,m^3-32\,m^2+79\,m-60)=0. \label{s5}
\end{equation}
Solving Eq. (\ref{s5}), we obtain
\begin{equation}
m \in \left \{ 0,\; {3\over 2},\;{5\over 2},\;4 \right \} .
\label{s6}
\end{equation}
Using the ansatz Eq. (\ref{s4}a) into Eq. (\ref{s3}) yields
\begin{equation}
y'' +{ A^2\,m\,(4-m)\over 20\, (A\,x+B)^2}\, y =0 . \label{s7}
\end{equation}
The general solution of Eq. (\ref{s7}) is
\begin{equation}
y = \left \{ \begin{array}{ll}
 k_1 + k_2\,(Ax+B) & \mbox{ if  } m\in \{ 0,4\},\\
 k_1\, (Ax+B)^{1/4}+ k_2 \,(Ax+B)^{3/4} & \mbox { if } m\in \{3/2
 ,5/2\},
 \end{array}  \right.
 \label{s8}
 \end{equation}
 where $k_1$ and $k_2$ are arbitrary constants.  We infer from Eq.
 (\ref{s2}) that

   \begin{equation}
g = \left \{ \begin{array}{ll}
 {L + M\,(Ax+B)\over P + Q\,(Ax+B)}  & \mbox{ if  } m\in \{ 0,4\},\\
 &  \\
 {L\, (Ax+B)^{1/4}+ M \,(Ax+B)^{3/4}\over P\, (Ax+B)^{1/4}+ Q \,(Ax+B)^{3/4}} & \mbox { if } m\in \{3/2
 ,5/2\},
 \end{array}  \right.
 \label{s9}
 \end{equation}
 where $L$, $M$, $P$ and $Q$ are constants satisfying $LQ-MP\ne 0$.

 The models defined by Eqs. (\ref{s4}a), (\ref{s6}) and (\ref{s9}) generalize
 Gottlieb's \cite{gott} iso-spectral models.

\subsection{The class $[ A_1\oplus A_2]$}
A particular solution of the system (\ref{d18})-(\ref{d19}) that
does not satisfy Eqs. (\ref{d11})-(\ref{d12}) is $f(x)=x,\;g(x)=x$.
Based on this particular solution, we choose as representative of $[
A_1\oplus A_2]$ the equation
\begin{equation}
(X\,U_{XX})_{XX}+X\,U_{TT}=0. \label{e8}
\end{equation}
The  finite portion of the Lie symmetry algebra of Eq. (\ref{e8}) is
generated by $Y_1$, $Y_2$, and $Y_{3,2}=4T\partial_T
+2X\partial_X-U\partial_U$. An invertible transformation (\ref{e2})
maps a element of $[ A_1\oplus A_2]$ to Eq. (\ref{e8}) if and only
if it maps $<X_1,X_2,X_3>$ to $<Y_1,Y_2,Y_{3,2}>$. We search for
such a transformation by mapping the basis elements as follows.
\begin{equation}
X_1\mapsto Y_1,\quad X_2\mapsto Y_2,\quad X_3\mapsto Y_{3,2}.
\label{e9}
\end{equation}
Simple calculations show that Eq. (\ref{e9}) is realized if and
only if
\begin{equation}
T=t+l_1\,g^2,\quad X=2\,l_2\,g,\quad U=l_3\,u\,\sqrt{f\,g'^3\over
g},\label{e10}
\end{equation}
where $l_1$, $l_2\ne 0$, and $l_3\ne 0$ are constants.

Thus we have proved the following theorem.
\begin{theo} Equations of $[A_2\oplus A_1]$ are homeomorphic
to Eq. (\ref{e8}). A transformation that maps an arbitrary element
of $[A_2\oplus A_1]$ to Eq. (\ref{e8}) is given by Eq.
(\ref{e10}).
\end{theo}
\subsection{The class $[3 A_1]$}
 A simple solution of Eqs. (\ref{d22})-(\ref{d23}) which does not satisfy Eqs. (\ref{d11})-(\ref{d12}) is $ f(x)=1,\;
 g(x)=\ln x$. To this solution corresponds the representative
 \begin{equation}
 U_{XXXX}+ X^{-4}\; U_{TT}=0. \label{e11}
 \end{equation}
The quotient symmetry Lie algebra  $S/L^{\infty}$ of Eq.
(\ref{e11}) is spanned by $Y_1$, $Y_2$, and
$Y_{3,2}=2\,X\,\partial_X+3\,U\,\partial_U$.  An element of $[3
A_1]$ will be homeomorphic to Eq. (\ref{e11}) provided  the Lie
algebras $<X_1,X_2,X_4>$ and $<Y_1,Y_2,Y_{3,2}>$ are similar. We
look for a similarity that transforms the basis vectors as
follows.
\begin{equation}
X_1 \mapsto Y_1,\quad X_2 \mapsto Y_1, \quad X_4 \mapsto Y_{3,2}.
\label{e12}
\end{equation}
Elementary reckoning shows that the mapping (\ref{e12}) is realized
by the transformation
\begin{equation}
T =t+m_1,\quad X=m_2\e^{g},\quad U= m_3u\sqrt{f g'^3\e^{3g}},
\label{e13}
\end{equation}
where $m_1$, $m_2\ne 0$, and $m_3\ne 0$ are constants.

 In summary, we have proved the following statement.
\begin{theo} Equations of $[3 A_1]$ are homeomorphic
to Eq. (\ref{e11}). A transformation that maps an arbitrary element
of $[3 A_1]$ to Eq. (\ref{e11}) is given by Eq. (\ref{e13}).
\end{theo}

\section{Conclusion}
 We have studied in details symmetry breaking of Euler-Bernoulli
 equation. We have shown that the Lie symmetry algebra of the
 Euler-Bernoulli equation is one of the following: $2A_1\oplus
 L^{\infty}$, $3A_1\oplus L^{\infty}$, $A_1\oplus A_2\oplus
 L^{\infty}$, or $A_{3,3}\oplus A_1\oplus L^{\infty}$, where
 $L^{\infty}$ is the infinite-dimensional Lie algebra spanned by
 solution symmetries. Equations admitting  a given symmetry class
 are characterized completely in terms of a determined set of
 non-linear ordinary differential equations that physical parameters
 (flexural rigidity and lineal mass
 density) must fulfill. Equations of the same class can be mapped to
 each other via invertible transformations. For each class we
 provided  a simple representative and we explicitly constructed
 similarity mappings. For the particular class of equations equivalent to the uniform
 Euler-Bernoulli equation, we rigorously  constructed  explicit non-trivial examples
 that extend and generalize Gottlieb's \cite{gott} iso-spectral
 models.

\newpage
\appendix
\section{Differential function appearing Eqs. (\ref{d4})-(\ref{d6})}
\label{ap1}
\begin{small}
\begin{eqnarray*}
H_{11} &= & \frac{6\,{f'}^2}{f} - \frac{6\,g\,{f'}^3}{{f}^2\,g'} -
   8\,f'' + \frac{10\,g\,f'\,f''}{f\,g'}
     -
   \frac{6\,g\,{f'}^2\,g''}{f\,{g'}^2}
   +
   \frac{8\,g\,f''\,g''}{{g'}^2} +
   \frac{30\,f\,{g''}^2}{{g'}^2} \\
   & & -
   \frac{60\,f\,g\,{g''}^3}{{g'}^4} -
   \frac{4\,g\,f^{(3)}}{g'}
    - \frac{20\,f\,g^{(3)}}{g'}
     +
   \frac{60\,f\,g\,g''\,g^{(3)}}{{g'}^3} -
   \frac{10\,f\,g\,g^{(4)}}{{g'}^2} \\
H_{12}&=& \frac{-6\,{f'  }^3}{{f  }^2\,g'  } +
   \frac{10\,f'  \,f''  }{f  \,g'  } -
   \frac{6\,{f'  }^2\,g''  }{f  \,{g'  }^2}
     +
   \frac{8\,f''  \,g''  }{{g'  }^2}
    -
   \frac{60\,f  \,{g''  }^3}{{g'  }^4}
    - \frac{4\,f^{(3)}  }{g'  }     +
   \frac{60\,f  \,g''  \,g^{(3)}  }{{g'  }^3}\\
   & &  -
   \frac{10\,f  \,g^{(4)}  }{{g'  }^2} \\
H_{21} & = & \frac{12\,{f'  }^4}{{f  }^3} -
   \frac{12\,g  \,{f'  }^5}{{f  }^4\,g'  } -
   \frac{28\,{f'  }^2\,f''  }{{f  }^2}
    +
   \frac{34\,g  \,{f'  }^3\,f''  }{{f  }^3\,g'  } +
   \frac{10\,{f''  }^2}{f  } -
   \frac{21\,g  \,f'  \,{f''  }^2}{{f  }^2\,g'  }
    \\
    & &  -
   \frac{12\,g  \,{f'  }^4\,g''  }{{f  }^3\,{g'  }^2}
    +
   \frac{6\,{f'  }^3\,g''  }{{f  }^2\,g'  } +
   \frac{28\,g  \,{f'  }^2\,f''  \,g''  }{{f  }^2\,{g'  }^2}
    -
   \frac{11\,f'  \,f''  \,g''  }{f  \,g'  } -
   \frac{10\,g  \,{f''  }^2\,g''  }{f  \,{g'  }^2} \\
   & & -
   \frac{12\,g  \,{f'  }^3\,{g''  }^2}{{f  }^2\,{g'  }^3}
    +
   \frac{6\,{f'  }^2\,{g''  }^2}{f  \,{g'  }^2}
    +
   \frac{22\,g  \,f'  \,f''  \,{g''  }^2}{f  \,{g'  }^3} -
   \frac{3\,f''  \,{g''  }^2}{{g'  }^2}
     -
   \frac{12\,g  \,{f'  }^2\,{g''  }^3}{f  \,{g'  }^4}\\
   & & -
   \frac{60\,f'  \,{g''  }^3}{{g'  }^3}  +
   \frac{6\,g  \,f''  \,{g''  }^3}{{g'  }^4}
   +
   \frac{120\,g  \,f'  \,{g''  }^4}{{g'  }^5} +
   \frac{180\,f  \,{g''  }^4}{{g'  }^4} -
   \frac{360\,f  \,g  \,{g''  }^5}{{g'  }^6}\\
   & &  +
   \frac{10\,f'  \,f^{(3)}  }{f  }
    -
   \frac{12\,g  \,{f'  }^2\,f^{(3)}  }{{f  }^2\,g'  }
   + \frac{9\,g  \,f''  \,f^{(3)}  }{f  \,g'  }
    -
   \frac{10\,g  \,f'  \,g''  \,f^{(3)}  }{f  \,{g'  }^2} +
   \frac{6\,g''  \,f^{(3)}  }{g'  } \\
   & & -
   \frac{12\,g  \,{g''  }^2\,f^{(3)}  }{{g'  }^3}
    +
   \frac{6\,g  \,{f'  }^3\,g^{(3)}  }{{f  }^2\,{g'  }^2} -
   \frac{4\,{f'  }^2\,g^{(3)}  }{f  \,g'  } -
   \frac{11\,g  \,f'  \,f''  \,g^{(3)}  }{f  \,{g'  }^2}
     +
   \frac{2\,f''  \,g^{(3)}  }{g'  } \\
   & & +
   \frac{12\,g  \,{f'  }^2\,g''  \,g^{(3)}  }{f  \,{g'  }^3} +
   \frac{70\,f'  \,g''  \,g^{(3)}  }{{g'  }^2}
    -
   \frac{6\,g  \,f''  \,g''  \,g^{(3)}  }{{g'  }^3} -
   \frac{180\,g  \,f'  \,{g''  }^2\,g^{(3)}  }{{g'  }^4}\\
    & & -
   \frac{300\,f  \,{g''  }^2\,g^{(3)}  }{{g'  }^3}
     +\frac{720\,f  \,g  \,{g''  }^3\,g^{(3)}  }{{g'  }^5} +
   \frac{6\,g  \,f^{(3)}  \,g^{(3)}  }{{g'  }^2}
    +
   \frac{30\,g  \,f'  \,{g^{(3)}  }^2}{{g'  }^3}\\
    & & + \frac{60\,f  \,{g^{(3)}  }^2}{{g'  }^2} -
   \frac{270\,f  \,g  \,g''  \,{g^{(3)}  }^2}{{g'  }^4}
     -
   4\,f^{(4)}   + \frac{3\,g  \,f'  \,f^{(4)}  }{f  \,g'  }
   +\frac{4\,g  \,g''  \,f^{(4)}  }{{g'  }^2} \\
   & & -
   \frac{2\,g  \,{f'  }^2\,g^{(4)}  }{f  \,{g'  }^2}
     -
   \frac{15\,f'  \,g^{(4)}  }{g'  }    +
   \frac{g  \,f''  \,g^{(4)}  }{{g'  }^2} +
   \frac{40\,g  \,f'  \,g''  \,g^{(4)}  }{{g'  }^3}
    +
   \frac{75\,f  \,g''  \,g^{(4)}  }{{g'  }^2}\\
    & & -   \frac{180\,f  \,g  \,{g''  }^2\,g^{(4)}  }{{g'  }^4}
    +\frac{60\,f  \,g  \,g^{(3)}  \,g^{(4)}  }{{g'  }^3}
    -
   \frac{g  \,f^{(5)}  }{g'  }
    -
   \frac{5\,g  \,f'  \,g^{(5)}  }{{g'  }^2}
   -   \frac{12\,f  \,g^{(5)}  }{g'  }\\
    & & +
      \frac{30\,f  \,g  \,g''  \,g^{(5)}  }{{g'  }^3}
     -
   \frac{3\,f  \,g  \,g^{(6)}  }{{g'  }^2}  \\
H_{22} &=& \frac{-12\,{f'  }^5}{{f  }^4\,g'  } +
   \frac{34\,{f'  }^3\,f''  }{{f  }^3\,g'  } -
   \frac{21\,f'  \,{f''  }^2}{{f  }^2\,g'  }
    -
   \frac{12\,{f'  }^4\,g''  }{{f  }^3\,{g'  }^2}
     +
   \frac{28\,{f'  }^2\,f''  \,g''  }{{f  }^2\,{g'  }^2} \\
   & & -
   \frac{10\,{f''  }^2\,g''  }{f  \,{g'  }^2}
    -
   \frac{12\,{f'  }^3\,{g''  }^2}{{f  }^2\,{g'  }^3} +
   \frac{22\,f'  \,f''  \,{g''  }^2}{f  \,{g'  }^3} -
   \frac{12\,{f'  }^2\,{g''  }^3}{f  \,{g'  }^4}
    +
   \frac{6\,f''  \,{g''  }^3}{{g'  }^4}\\
   & &  +
   \frac{120\,f'  \,{g''  }^4}{{g'  }^5} -
   \frac{360\,f  \,{g''  }^5}{{g'  }^6}
    -
   \frac{12\,{f'  }^2\,f^{(3)}  }{{f  }^2\,g'  } +
   \frac{9\,f''  \,f^{(3)}  }{f  \,g'  }
    -
   \frac{10\,f'  \,g''  \,f^{(3)}  }{f  \,{g'  }^2}\\
   & &  -
   \frac{12\,{g''  }^2\,f^{(3)}  }{{g'  }^3} +
   \frac{6\,{f'  }^3\,g^{(3)}  }{{f  }^2\,{g'  }^2} -
   \frac{11\,f'  \,f''  \,g^{(3)}  }{f  \,{g'  }^2}
    +
   \frac{12\,{f'  }^2\,g''  \,g^{(3)}  }{f  \,{g'  }^3}
     -
   \frac{6\,f''  \,g''  \,g^{(3)}  }{{g'  }^3}\\
    & & -
   \frac{180\,f'  \,{g''  }^2\,g^{(3)}  }{{g'  }^4}
    +
   \frac{720\,f  \,{g''  }^3\,g^{(3)}  }{{g'  }^5}
   +
   \frac{6\,f^{(3)}  \,g^{(3)}  }{{g'  }^2} +
   \frac{30\,f'  \,{g^{(3)}  }^2}{{g'  }^3}\\
   & &  -
   \frac{270\,f  \,g''  \,{g^{(3)}  }^2}{{g'  }^4}
    +
   \frac{3\,f'  \,f^{(4)}  }{f  \,g'  } +
   \frac{4\,g''  \,f^{(4)}  }{{g'  }^2}
    -
   \frac{2\,{f'  }^2\,g^{(4)}  }{f  \,{g'  }^2} +
   \frac{f''  \,g^{(4)}  }{{g'  }^2}
    +
   \frac{40\,f'  \,g''  \,g^{(4)}  }{{g'  }^3}\\
   & & -
   \frac{180\,f  \,{g''  }^2\,g^{(4)}  }{{g'  }^4}
    +
   \frac{60\,f  \,g^{(3)}  \,g^{(4)}  }{{g'  }^3} -
   \frac{f^{(5)}  }{g'  } - \frac{5\,f'  \,g^{(5)}  }{{g'  }^2}
   +
   \frac{30\,f  \,g''  \,g^{(5)}  }{{g'  }^3} -
   \frac{3\,f  \,g^{(6)}  }{{g'  }^2}
   \end{eqnarray*}
   \end{small}
   \section{Successive reduction of the order of Eq. (\ref{d12})} \label{ap2}
 The Lie brackets of the symmetry operators are
 $[\Gamma_1,\Gamma_2]=
 \Gamma_1,\;[\Gamma_1,\Gamma_3]=0,\;[\Gamma_2,\Gamma_3]=0$. Thus
 $<\Gamma_1>$ and $<\Gamma_3>$ are ideals of the symmetry Lie
 algebra. We may start reduction using $\Gamma_1$ or $\Gamma_3$. It
 is crucial start reduction using an ideal of the symmetry Lie
 algebra in order to preserve the remaining symmetries.

A basis of first-order differential invariants of $\Gamma_1$ is
formed by $ f $ and $f'$. We define new dependent and independent
variables by
\begin{equation}
y=f',\quad t=f. \label{ar1}
\end{equation}
In the new variables, Eq. (\ref{d12}) reads
\begin{equation}
y^{(3)} ={\ddot y\over t} -4\,{\dot y\, \ddot y\over y}-{12\over
5}\, {\dot y \over t^2}+{21\over 20}\,{\dot y^2\over y}-7\;{\dot
y^3\over y^2} +{9\over 10}\,{y\over t^3} , \label{ar2}
\end{equation}
where the over-dot stands for differentiation with respect to $t$.
In the new variables, the symmetries $\Gamma_2$ and $\Gamma_3$ are
$$ \Gamma_2 =-y\,\partial_y\,,\quad \Gamma_3 =t\,\partial_t -y\,\partial_y
\,.
$$
It can be verify that Eq. (\ref{ar2}) inherits the symmetries
$\Gamma_2$ and $\Gamma_3$ as expected.

A basis of first-order differential invariant of $\Gamma_2$ is
formed by $t$ and
\begin{equation}
z ={\dot y \over y}\, \cdot \label{ar3}
\end{equation}
In the new variables $t$ and $z$, Eq. (\ref{ar2}) becomes
\begin{equation}
\ddot z ={\dot z\over t}-7\,z\,\dot z -{12\over 5}\,z+{41\over
20}\,{z^2\over t}-16\,z^3+{9\over 10}\,{1\over t^3}\, \cdot
\label{ar4}
\end{equation}
In  the variables $t$ and $z$, the symmetry operator $\Gamma_3$
becomes
$$ \Gamma_3 =t\partial_t-z\,\partial_z \,.$$
A basis of first-order differential invariant of $\Gamma_3$ is
\begin{equation}
u=t\,z,\quad v=tz+t^2\,z' . \label{ar4}
\end{equation}

In the new variables $u$ and $v$, Eq. (\ref{ar4}) is

\begin{equation}
{dv\over du} =5-7u -{320\, u^3-181\,u^2+108\,u-18 \over 20\, v} \,
\cdot \label{ar5}
\end{equation}
Eq. (\ref{ar5}) is an Abelian equation of the second kind.


\begin{thebibliography}{99}
\bibitem{gott} H. P. W. Gottlieb, Isospectral Euler-Bernoulli beam with continuous density and rigidity
 functions, Proceedings of the Royal Society London A 413 (1987) 235-250.
\bibitem{blku} G. Bluman and S. Kumei, Symmetries and Differential Equations, Springer-Verlag, New York, 1989.
 \bibitem{ovsi} L.V. Ovsiannikov, Group analysis of differential
equations, Academic Press, New York, 1982.
\bibitem{olve} P. Olver,  Applications of Lie groups to Differential Equations, Springer-Verlag, New
York, 1993.
\bibitem{pawi} J. Patera and P. Winternitz, Subalgebras of real three- and four-dimensional Lie
algebras, Journal of  Mathematical  Physics 18 (1977) 1449-1455.
\bibitem{inc} E. L. Ince, Differential equations, Dover, New York,
1956.
\end{thebibliography}
\end{document}